\documentclass[10pt]{article}
\usepackage{lmodern}
\usepackage{amsmath}
\usepackage[T1]{fontenc}
\usepackage[utf8]{inputenc}
\usepackage{authblk}
\usepackage{amsfonts}
\usepackage{graphicx}
\usepackage{tkz-euclide}
\tikzset{elegant/.style={smooth,thick,samples=50,cyan}}
\tikzset{eaxis/.style={->,>=stealth}}
\usepackage{color,tikz}
\usetikzlibrary{decorations.pathreplacing,calc}
\usepackage{rotating}
\usepackage{amssymb}
\usepackage[english]{babel}
\usepackage{color}
\usepackage{amsthm}
\usepackage{graphicx}
\usepackage{mathrsfs}
\usepackage{makecell}
\usepackage{microtype}
\usepackage{mathscinet}
\usepackage{array}
\usepackage{multirow}
\usepackage{enumerate}
\usepackage[cal=boondoxo,bb=ams]{mathalfa}
\usepackage{hyperref}
\hypersetup{hidelinks}
\newtheorem{defn}{Definition}[section]
\newtheorem{theorem}{Theorem}[section]

\newtheorem{remark}{Remark}[section]

\newcommand{\ml}{\mathcal}
\newcommand{\mb}{\mathbb}

\title{Blow-up and lifespan estimates for Nakao's type problem with nonlinearities of derivative type}
\author[1]{Wenhui Chen\thanks{Corresponding author: Wenhui Chen (wenhui.chen.math@gmail.com)}}
\affil[1]{School of Mathematical Sciences, Shanghai Jiao Tong University, 200240 Shanghai, China}
\date{}
\begin{document}
		\maketitle
\begin{abstract}
In the present paper, we investigate blow-up and lifespan estimates for a class of semilinear hyperbolic coupled system in $\mathbb{R}^n$ with $n\geqslant 1$, which is part of the so-called Nakao's type problem weakly coupled a semilinear damped wave equation with a semilinear wave equation with nonlinearities of derivative type. By constructing two time-dependent functionals and employing an iteration method for unbounded multiplier with slicing procedure, the results of blow-up and upper bound estimates for the lifespan of energy solutions are derived. The model seems to be hyperbolic-like instead of parabolic-like. Particularly, the blow-up result for one dimensional case is optimal.
	\medskip\\
	\textbf{AMS Classification (2010)} Primary 35L52; Secondary 35B44
\end{abstract}
\fontsize{12}{15}
\selectfont

%\tableofcontents
\section{Introduction} \label{Sec1}
 The problem of critical curve, which describes the threshold condition between global (in time) existence of small data weak solutions and blow-up of small data weak solutions, of the power exponents for the weakly coupled system of wave equations and damped wave equations was proposed by Professor Mitsuhiro Nakao, Emeritus of Kyushu University (see, for instance, \cite{NishiharaWakasugi2015,Wakasugi2017}), namely,
\begin{equation}\label{General Nakao's type problem}
\begin{cases}
u_{tt}-\Delta u+u_t=f_1(v,v_t),&x\in\mb{R}^n, \ t>0,\\
v_{tt}-\Delta v=f_2(u,u_t),&x\in\mb{R}^n,\ t>0,\\
(u,u_t,v,v_t)(0,x)=(u_0,u_1,v_0,v_1)(x),&x\in\mb{R}^n,
\end{cases}
\end{equation}
where the nonlinearties on the right-hand sides are given by the mixture of the power type and the derivative type nonlinearities
\begin{align*}
f_1(v,v_t)&:=d_1|v|^{p_1}+d_2|v_t|^{p_2},\\
f_2(u,u_t)&:=d_3|u|^{q_1}+d_4|u_t|^{q_2},
\end{align*}
carrying some nonnegative constants $d_1,\dots ,d_4$ and $p_1,p_2,q_1,q_2>1$. Here, to guarantee  the hyperbolic coupled system \eqref{General Nakao's type problem} being a nonlinear problem, we have to restrict ourselves that the coefficients satisfy $d_1+d_2\neq0$ and $d_3+d_4\neq0$. Roughly speaking, the main difficulty to treat Nakao's type problem is to understand varying degrees of influence from damped wave equations and wave equations. It is well-known that $L^2$ decay properties and diffusion phenomenon hold in damped wave equations due to the frictional damping $u_t$. However, these effects disappear in wave equations and Huygens' principle and energy conservation are valid in  wave equations, which make huge differences of the treatments between semilinear wave equations and semilinear damped wave equations. In other words, Nakao's type problem bridges a connection between semilinear damped wave equations and semilinear wave equations through weakly coupled form in the power nonlinearities. We emphasize that the critical condition for the Cauchy problem \eqref{General Nakao's type problem} is still an open problem for any $d_1+d_2\neq0$ and $d_3+d_4\neq0$.

In recent years, some blow-up results for Nakao's type problem with power nonlinearities, namely,
\begin{equation}\label{Nakao's problem power}
\begin{cases}
u_{tt}-\Delta u+u_t=|v|^p,&x\in\mb{R}^n, \ t>0,\\
v_{tt}-\Delta v=|u|^q,&x\in\mb{R}^n,\ t>0,\\
(u,u_t,v,v_t)(0,x)=(u_0,u_1,v_0,v_1)(x),&x\in\mb{R}^n,
\end{cases}
\end{equation}
that is the special case of the hyperbolic coupled system \eqref{General Nakao's type problem} with $d_1=d_3=1$, $d_2=d_4=0$ and $p_1=p$, $q_1=q$, have been derived in \cite{Wakasugi2017,Chen-Reissig-2020}. With the aim of guaranteeing existence of local (in time) solutions, let us generally assume $p,q\leqslant n/(n-2)$ if $n\geqslant 3$ in the discussion of this paragraph. Firstly, by using a test function method, the author of \cite{Wakasugi2017} demonstrated blow-up of local (in time) weak solutions with suitable assumption on initial data providing that
\begin{align}\label{Condition Wakasugi}
\max\left\{\frac{q/2+1}{pq-1}+\frac{1}{2},\frac{q+1}{pq-1},\frac{p+1}{pq-1}\right\}\geqslant\frac{n}{2}.
\end{align}
The condition \eqref{Condition Wakasugi} is optimal in $n=1$ since it is equivalent to $1<p,q<\infty$. Later, the authors of \cite{Chen-Reissig-2020} employed an iteration method associated with slicing procedure to improve the blow-up condition \eqref{Condition Wakasugi} for $n\geqslant 2$. Precisely, the result in \cite{Wakasugi2017} was partially improved for $n=2,3$ and completely improved for $n\geqslant 4$ such that if 
\begin{align}\label{Condition Chen-Reissig}
\alpha_{\mathrm{N}}(p,q):=\max\left\{\frac{q/2+1}{pq-1},\frac{2+p^{-1}}{pq-1},\frac{1/2+p}{pq-1}-\frac{1}{2}\right\}>\frac{n-1}{2},
\end{align}
then every non-zero local (in time) energy solution blows up in finite time. In other words, the authors of \cite{Chen-Reissig-2020} observed that Nakao's type problem with power nonlinearities is hyperbolic-like model rather than parabolic-like model due to fact that the component
\begin{align*}
	\frac{2+p^{-1}}{pq-1}>\frac{n-1}{2}
\end{align*}
 plays a dominant role when $n\geqslant 3$. This effect comes from the semilinear wave equations. More detail explanations of parabolic-like versus hyperbolic-like are referred  interested readers to Section 2.1 in \cite{Chen-Reissig-2020}. Therefore, an interesting and viable problem is to ask the situation of nonlinearities of derivative type, i.e. the hyperbolic coupled system \eqref{General Nakao's type problem} with $d_1=d_3=0$, $d_2=d_4=1$ and $p_2=p$, $q_2=q$. At this time, one may notice that the time-derivative of solution not only exists in the linear part (the damped wave equation), but also appears in the nonlinear parts of both equations. We would like to understand:
 \begin{center}
 	Do the nonlinear terms including the time-derivative of solutions change the model from hyperbolic-like (i.e. the nonlinear problem (\ref{Nakao's problem power})) to parabolic-like?
 \end{center}
  We will give a possible answer from the blow-up point of view  that Nakao's type problem with nonlinearities of derivative type still could be hyperbolic-like model.

In this paper, we study blow-up of solutions and lifespan estimates from the above for Nakao's type problem with derivative type nonlinearities, namely,
\begin{equation}\label{Nakao's problem derivative}
\begin{cases}
u_{tt}-\Delta u+u_t=|v_t|^p,&x\in\mb{R}^n, \ t>0,\\
v_{tt}-\Delta v=|u_t|^q,&x\in\mb{R}^n,\ t>0,\\
(u,u_t,v,v_t)(0,x)=\varepsilon(u_0,u_1,v_0,v_1)(x),&x\in\mb{R}^n,
\end{cases}
\end{equation}
where $p,q>1$ and $\varepsilon$ is a positive parameter describing the size of initial data. As we will show in Theorem \ref{Thm 1}, the blow-up condition of Nakao's type problem \eqref{Nakao's problem derivative} is strongly related to the Glassey exponent, which is the critical exponent for the semilinear wave equations with derivative type nonlinearity (see \cite{HidanoWangYokoyama2012} for some detail introductions). The approach to derive our result is mainly based on an iteration method for unbounded multiplier with  slicing procedure by setting suitable time-dependent functionals. 

Let us show some results for semilinear wave equations and semilinear damped wave equations, which are strongly related to Nakao's type problem \eqref{Nakao's problem derivative}. Concerning the weakly coupled system of semilinear wave equations
\begin{equation}\label{Wave coupled derivative}
\begin{cases}
u_{tt}-\Delta u=|v_t|^p,&x\in\mb{R}^n, \ t>0,\\
v_{tt}-\Delta v=|u_t|^q,&x\in\mb{R}^n,\ t>0,\\
(u,u_t,v,v_t)(0,x)=(u_0,u_1,v_0,v_1)(x),&x\in\mb{R}^n,
\end{cases}
\end{equation}
the critical curve is given by
\begin{align}\label{Critical Curve Wave derivative}
\alpha_{\mathrm{W}}(p,q):=\frac{\max\{p,q\}+1}{pq-1}=\frac{n-1}{2}.
\end{align}
Particularly, under certain integral sign assumptions for initial data, if $\alpha_{\mathrm{W}}(p,q)\geqslant (n-1)/2$, then every non-trivial local (in time) solution $(u,v)$ blows up in finite time. Considering the critical curve \eqref{Critical Curve Wave derivative}, we refer to the related works \cite{Deng-1999,Xu,Kubo-Kubota-Sunagawa,Ikeda-Sobajima-Wakasa,PalTak19dt}. Taking our  consideration of the special case $p=q$, the critical exponent is given by the so-called Glassey exponent
\begin{align*}
p_{\mathrm{Gla}}(n):=\begin{cases}
\infty&\mbox{if}\ \ n=1,\\
\displaystyle{\frac{n+1}{n-1}}&\mbox{if}\ \ n\geqslant 2,
\end{cases}
\end{align*}
which is also the critical exponent for the single semilinear wave equation with nonlinearity $|u_t|^p$. One may see the  validity of the Glassey exponent in \cite{John1981,Sideris1983,Masuda1983,Schaeffer1986,Rammaha1987,Agemi1991,HidanoTsutaya1995,Tzvetkov1998,Zhou2001,HidanoWangYokoyama2012,Lai-Takamura-2019} and reference therein. Next, we turn to the weakly coupled system of semilinear classical damped wave equations with nonlinearities of derivative type as follows:
\begin{equation}\label{Damped Wave coupled derivative}
\begin{cases}
u_{tt}-\Delta u+u_t=|v_t|^p,&x\in\mb{R}^n, \ t>0,\\
v_{tt}-\Delta v+v_t=|u_t|^q,&x\in\mb{R}^n,\ t>0,\\
(u,u_t,v,v_t)(0,x)=(u_0,u_1,v_0,v_1)(x),&x\in\mb{R}^n.
\end{cases}
\end{equation}
According to the previous study \cite{D'Ab-Ebert-2017}, we believe that the global (in time) solution of the last system uniquely exists  for any $1<p,q<\infty$ with $p,q\leqslant n/(n-2)$ if $n\geqslant 3$. Namely, the solution does not blow up for any dimensions. This effect also appears in the wave equations with scale-invariant damping of the effective case (see, among other things in Theorem 2.2 of \cite{Palmieri-Tu-2019} by letting parameters $\mu_1,\mu_2\to\infty$). Thus, the consideration of Nakao's type problem \eqref{Nakao's problem derivative} is reasonable.

\medskip

\noindent\textbf{Notation: } We give some notations to be used in this paper. We write $f\lesssim g$ when there exists a positive constant $C$ such that $f\leqslant Cg$. We denote $\lceil r\rceil:=\min\{C\in\mb{Z}: r\leqslant C\}$ as the ceiling function. Moreover, $B_R$ denotes the ball around the origin with radius $R$ in $\mathbb{R}^n$.

\section{Main result}
Let us first introduce a suitable definition of energy solutions of Nakao's type problem \eqref{Nakao's problem derivative}.
\begin{defn}\label{Defn Energy Solution}
	Let $(u_0,u_1,v_0,v_1)\in \left(H^1(\mb{R}^n)\times L^2(\mb{R}^n)\right)\times \left(H^1(\mb{R}^n)\times L^2(\mb{R}^n)\right)$. One may say that $(u,v)$ is an energy solution of Nakao's type problem \eqref{Nakao's problem derivative} on $[0,T)$ if
	\begin{align*}
	u\in\ml{C}([0,T),H^1(\mb{R}^n))\cap \ml{C}^1([0,T),L^2(\mb{R}^n))\ \ \mbox{and}\ \ u_t\in L_{\mathrm{loc}}^q([0,T)\times\mb{R}^n),\\
	v\in\ml{C}([0,T),H^1(\mb{R}^n))\cap \ml{C}^1([0,T),L^2(\mb{R}^n))\ \ \mbox{and}\ \ v_t\in L_{\mathrm{loc}}^p([0,T)\times\mb{R}^n),
	\end{align*}
	satisfies $(u,v)(0,\cdot)=\varepsilon(u_0,v_0)$ in $H^1(\mb{R}^n)\times H^1(\mb{R}^n)$ and the following integral relations:
	\begin{align}\label{Defn.Energy.Solution 01}
	&\int_0^t\int_{\mb{R}^n}(-u_t(s,x)\phi_s(s,x)+u_t(s,x)\phi(s,x)+\nabla u(s,x)\cdot\nabla\phi(s,x))\mathrm{d}x\mathrm{d}s\notag\\
	&+\int_{\mb{R}^n}u_t(t,x)\phi(t,x)\mathrm{d}x-\varepsilon\int_{\mb{R}^n}u_1(x)\phi(0,x)\mathrm{d}x=\int_0^t\int_{\mb{R}^n}|v_t(s,x)|^p\phi(s,x)\mathrm{d}x\mathrm{d}s
	\end{align}
	and
	\begin{align}\label{Defn.Energy.Solution 02}
	&\int_0^t\int_{\mb{R}^n}(-v_t(s,x)\psi_s(s,x)+\nabla v(s,x)\cdot\nabla\psi(s,x))\mathrm{d}x\mathrm{d}s\notag\\
	&+\int_{\mb{R}^n}v_t(t,x)\psi(t,x)\mathrm{d}x-\varepsilon\int_{\mb{R}^n}v_1(x)\psi(0,x)\mathrm{d}x=\int_0^t\int_{\mb{R}^n}|u_t(s,x)|^q\psi(s,x)\mathrm{d}x\mathrm{d}s
	\end{align}
hold	for any test functions $\phi,\psi\in\ml{C}_0^{\infty}([0,T)\times\mb{R}^n)$ and any $t\in(0,T)$.
\end{defn}
\begin{remark}
Actually, similarly to treatments in \cite{Chen-Reissig-2020}, by applying further steps of integration by parts in \eqref{Defn.Energy.Solution 01} as well as \eqref{Defn.Energy.Solution 02}, respectively, and taking $t\to T$, we may claim that $(u,v)$ fulfills the definition of weak solutions of Nakao's type problem \eqref{Nakao's problem derivative}.
\end{remark}
%Therefore, further steps of integration by parts in \eqref{Defn.Energy.Solution 01} and \eqref{Defn.Energy.Solution 02}, respectively, lead to
%\begin{align*}
%&\int_0^t\int_{\mb{R}^n}u(s,x)(\phi_{ss}(s,x)-\Delta\phi(s,x)-\phi_s(s,x))\mathrm{d}x\mathrm{d}s\\
%&+\int_{\mb{R}^n}(u_t(t,x)\phi(t,x)+u(t,x)\phi(t,x)-u(t,x)\phi_s(t,x))\mathrm{d}x\\
%&-\int_{\mb{R}^n}(u_1(x)\phi(0,x)+u_0(x)\phi(0,x)-u_0(x)\phi_s(0,x))\mathrm{d}x=\int_0^t\int_{\mb{R}^n}|v_t(s,x)|^p\phi(s,x)\mathrm{d}x\mathrm{d}s
%\end{align*}
%and
%\begin{align*}
%&\int_0^t\int_{\mb{R}^n}v(s,x)(\psi_{ss}(s,x)-\Delta\psi(s,x))\mathrm{d}x\mathrm{d}s+\int_{\mb{R}^n}(v_t(t,x)\psi(t,x)-v(t,x)\psi_s(t,x))\mathrm{d}x\\
%&-\int_{\mb{R}^n}(v_1(x)\psi(0,x)-v_0(x)\psi_s(0,x))\mathrm{d}x=\int_0^t\int_{\mb{R}^n}|u_t(s,x)|^q\psi(s,x)\mathrm{d}x\mathrm{d}s.
%\end{align*}
%Letting $t\to T$, we claim that $(u,v)$ fulfills the definition of weak solutions to the hyperbolic coupled system \eqref{Nakao's problem derivative}.

From Banach's fixed point theorem and Duhamel's principle associated with some $L^2$ estimates of solutions of the corresponding linear Cauchy problem to \eqref{Nakao's problem derivative}, one may derive local (in time) existence of weak solutions with compact support localized in a ball with radius $R+t$ of Nakao's type problem \eqref{Nakao's problem derivative} with compactly supported data in a ball with radius $R$ if $p,q>1$ for $n=1,2$, and $1<p,q\leqslant n/(n-2)$ for $n\geqslant 3$.

Let us state the blow-up result for Nakao's type problem \eqref{Nakao's problem derivative}.
\begin{theorem}\label{Thm 1}
	Let us consider the exponents $p,q>1$ such that
	\begin{align}\label{Assumption.Condition}
	pq<\begin{cases}
	\infty&\mbox{if}\ \ n=1,\\
	\displaystyle{p_{\mathrm{Gla}}(n)}&\mbox{if}\ \ n\geqslant2.
	\end{cases}
	\end{align}
	Furthermore, let $(u_0,u_1,v_0,v_1)\in (H^1(\mb{R}^n)\times L^2(\mb{R}^n))\times (H^1(\mb{R}^n)\times L^2(\mb{R}^n))$ are nonnegative and compactly supported functions with supports contained in $B_R$ for some $R>0$ such that $u_1,v_1$ are not identically zero. Let $(u,v)$ be the local (in time) energy solution of Nakao's type problem \eqref{Nakao's problem derivative} according to Definition \ref{Defn Energy Solution} with lifespans $T=T(\varepsilon)$. Then, these solutions satisfy
	\begin{align}\label{Supp_uv}
	\mathrm{supp}\, u,\,\mathrm{supp}\,v\subset \{(t,x)\in[0,T)\times\mb{R}^n:|x|\leqslant R+t\}.
	\end{align}
	Moreover, there exists a positive constant $\varepsilon_0=\varepsilon_0(u_0,u_1,v_0,v_1,p,q,n,R)$ such that for any $\varepsilon\in(0,\varepsilon_0]$ the energy solution $(u,v)$ blows up in finite time. In addition, the upper bound estimate for the lifespans
	\begin{align*}
	T(\varepsilon)\leqslant C\varepsilon^{-\frac{2(pq-1)}{(n+1)-(n-1)pq}}
	\end{align*}
	holds, where $C>0$ is a constant independent of $\varepsilon$.
\end{theorem}
\begin{remark}
Due to the blow-up condition $pq<(n+1)/(n-1)$ for $n\geqslant2$, it is trivial that $p,q\leqslant n/(n-2)$ for $n\geqslant3$. In other words, under the condition of exponents $p$ and $q$ in Theorem \ref{Thm 1}, the weak solution having compact support in $B_{R+t}$ of Nakao's type problem \eqref{Nakao's problem derivative} locally (in time) exists.
\end{remark}
\begin{remark}
In the one dimensional case, the energy solution of Nakao's type problem \eqref{Nakao's problem derivative} blows up for all $1<p,q<\infty$, which implies the optimality for $n=1$. However, for the high dimensional cases when $n\geqslant2$, the critical curve in the $p-q$ plane for Nakao's type problem \eqref{Nakao's problem derivative}  is open.
\end{remark}

\begin{remark}
Due to the fact that
\begin{align*}
	\left\{(p,q):pq<\frac{n+1}{n-1}\right\}=\left\{(p,q):\frac{1}{pq-1}>\frac{n-1}{2}\right\} \subset \left\{(p,q):\alpha_{\mathrm{N}}(p,q)>\frac{n-1}{2}\right\}
\end{align*}
for any $n\geqslant 2$,  we may claim that the nonlinear terms including time-derivative of solutions, i.e. $(|v_t|^p,|u_t|^q)^{\mathrm{T}}$, will weaken the blow-up range of exponents for the Nakao's type model in the $p-q$ plane. However, the Nakao's type problem with nonlinearities of derivative type still behaves as a hyperbolic-like model whose reason will be shown later. For this reason, combing the explanation of Theorem \ref{Thm 1} and those in the recent paper \cite{Chen-Reissig-2020}, we may conjecture that for the general Nakao's type problem \eqref{General Nakao's type problem} with any $d_1+d_2\neq0$ and $d_3+d_4\neq0$, the model still has the behavior of hyperbolic-like.
\end{remark}

\subsection{Hyperbolic-like versus parabolic-like}
In this part, let us give some remarks and explanations on the blow-up conditions for energy solutions of Nakao's type problem \eqref{Nakao's problem derivative} with respect to the exponents $p,q$ for $n=2$ and $n\geqslant 3$, respectively. In the forthcoming discussion, we assume $p,q>1$.
\begin{figure}[http]
	\centering
	\begin{tikzpicture}
	\draw[->] (-0.2,0) -- (5.8,0) node[below] {$p$};
	\draw[->] (0,-0.2) -- (0,5.4) node[left] {$q$};
	\node[left] at (0,-0.2) {{$0$}};
	\node[below] at (1,0) {{$1$}};
	\node[left] at (0,1) {{$1$}};
	\node[left] at (0,3) {{$3$}};
	\node[below] at (3,0) {{$3$}};
	\node[left, color=red] at (5.9,3.8) {{$\longleftarrow$ $\frac{\max\{p,q\}+1}{pq-1}=\frac{1}{2}$}};
	\node[left, color=blue] at (3.3,2.2) {{$\longleftarrow$ $pq=3$}};
	\draw[dashed, color=black]  (0, 1)--(5.6, 1);
	\draw[dashed, color=black]  (1, 0)--(1, 5.4);
	\draw[dashed, color=black]  (3, 0)--(3, 3);
	\draw[dashed, color=black]  (0, 3)--(3, 3);
	\draw[color=blue] plot[smooth, tension=.7] coordinates {(1,3) (1.732,1.732) (3,1)};
	\draw[color=red] plot[smooth, tension=.7] coordinates {(2.6,5) (3,3) (5,2.6)};
	\node[left] at (4,-1) {{Case $n=2$}};
	%%%%%%%
	\draw[->] (7.8,0) -- (13.8,0) node[below] {$p$};
	\draw[->] (8,-0.2) -- (8,5.4) node[left] {$q$};
	\node[left] at (8,-0.2) {{$0$}};
	\node[below] at (9,0) {{$1$}};
	\node[left] at (8,1) {{$1$}};
	\node[below] at (11.5,0) {{$\frac{n}{n-2}$}};
	\node[left] at (8,3.5) {{$\frac{n}{n-2}$}};
		\node[below] at (10.5,0) {{$\frac{n+1}{n-1}$}};
	\node[left] at (8,2.5) {{$\frac{n+1}{n-1}$}};
	\node[left, color=red] at (13.8,3) {{$\longleftarrow$ $\frac{\max\{p,q\}+1}{pq-1}=\frac{n-1}{2}$}};
	\node[left, color=blue] at (12.1,1.5) {{$\longleftarrow$ $pq=\frac{n+1}{n-1}$}};
	\draw[dashed, color=black]  (8, 1)--(13.6, 1);
	\draw[dashed, color=black]  (9, 0)--(9, 5.4);
	\draw[dashed, color=black]  (8, 3.5)--(13.6, 3.5);
	\draw[dashed, color=black]  (11.5, 0)--(11.5, 5.4);
	\draw[dashed, color=black]  (10.5, 0)--(10.5, 2.5);
		\draw[dashed, color=black]  (8, 2.5)--(10.5, 2.5);
	\draw[color=blue] plot[smooth, tension=.7] coordinates {(9,2.5) (9.575,1.575) (10.5,1)};
	\draw[color=red] plot[smooth, tension=.7] coordinates {(10,3.5) (10.5,2.5) (11.5,2)};
	\node[left] at (12,-1) {{Case $n\geqslant 3$}};
	\end{tikzpicture}
	\caption{Blow-up conditions in the $p-q$ plane}
	\label{imggg}
\end{figure}
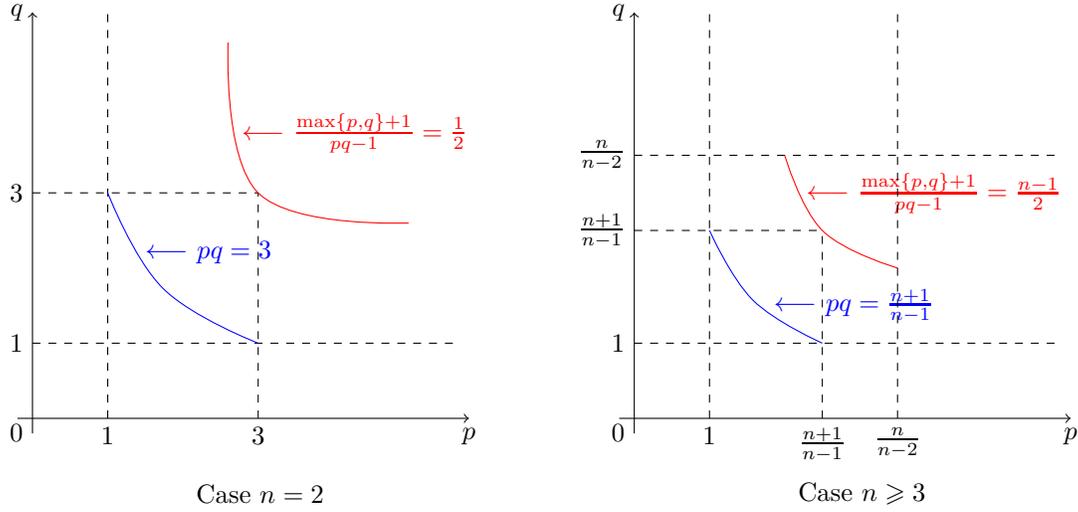
 
According to Figure \ref{imggg}, we may observe that
\begin{align}\label{Rela_Cond}
\left\{(p,q):pq<p_{\mathrm{Gla}}(n)=\frac{n+1}{n-1}\right\}\subset\left\{(p,q):\alpha_{\mathrm{W}}(p,q)=\frac{\max\{p,q\}+1}{pq-1}>\frac{n-1}{2}\right\}
\end{align}
for any $n\geqslant 2$. Again, $\alpha_{\mathrm{W}}(p,q)=(n-1)/2$ is the critical curve in the $p-q$ plane for the weak coupled system \eqref{Wave coupled derivative}. The effect \eqref{Rela_Cond} is caused by the influence of friction $u_t$ on the first equation of the Cauchy problem \eqref{Nakao's problem derivative}. For the reason of the blow-up condition $pq<p_{\mathrm{Gla}}(n)$, where $p_{\mathrm{Gla}}(n)$ is the critical exponent for semilinear wave equation with derivative type nonlinearity, we feel that Nakao's type problem \eqref{Nakao's problem derivative} is of hyperbolic-like (wave behavior) rather than parabolic-like. 

%Furthermore, we may rewrite the condition \eqref{Assumption.Condition} by
%\begin{align*}
%	\frac{1}{pq-1}>\frac{n-1}{2},
%\end{align*}
%which has the similar form to blow-up for the coupled system \eqref{Wave coupled derivative} in the subcritical case (see the right-hand side of \eqref{Rela_Cond}).
%
%Let us define a weak global (in time) solution to \eqref{Nakao's problem derivative}.
%\begin{defn}
%\begin{align}\label{Weak Solution 1}
%&-\int_0^{\infty}\int_{\mb{R}^n}u_t(t,x)\left(\phi_t(t,x)+\int_t^{\infty}\Delta\phi(\tau,x)\mathrm{d}\tau-\phi(t,x)\right)\mathrm{d}x\mathrm{d}t\notag\\
%&=\int_{\mb{R}^n}u_0(x)\int_0^{\infty}\phi(\tau,x)\mathrm{d}\tau\mathrm{d}x+\int_{\mb{R}^n}u_1(x)\phi(0,x)\mathrm{d}x+\int_0^{\infty}\int_{\mb{R}^n}|v_t(t,x)|^p\phi(t,x)\mathrm{d}x\mathrm{d}t
%\end{align}
%and
%\begin{align}\label{Weak Solution 2}
%&-\int_0^{\infty}\int_{\mb{R}^n}v_t(t,x)\left(\psi_t(t,x)+\int_t^{\infty}\Delta\psi(\tau,x)\mathrm{d}\tau\right)\mathrm{d}x\mathrm{d}t\notag\\
%&=\int_{\mb{R}^n}v_0(x)\int_0^{\infty}\psi(\tau,x)\mathrm{d}\tau\mathrm{d}x+\int_{\mb{R}^n}v_1(x)\psi(0,x)\mathrm{d}x+\int_0^{\infty}\int_{\mb{R}^n}|u_t(t,x)|^q\psi(t,x)\mathrm{d}x\mathrm{d}t
%\end{align}
%\end{defn}

\subsection{A shrift magnitude of the curve for blow-up conditions}
We observe an interesting effect in the blow-up conditions for Nakao's type problem, namely, the shrift curve with a certain magnitude from the weakly coupled system of wave equations with corresponding nonlinearities. Later, we assume $p,q\leqslant n/(n-2)$ if $n\geqslant 3$ and $p,q>1$ for any $n\geqslant 1$ to guarantee local (in time) existence of solutions.

Let us first consider the blow-up result for Nakao's type problem \eqref{Nakao's problem power} with power nonlinearities, i.e. $d_1=d_3=1$ and $d_2=d_4=0$ in the Cauchy problem \eqref{General Nakao's type problem}. The authors in \cite{Chen-Reissig-2020} proved blow-up of energy solutions for $n\geqslant 3$ if the exponents $p,q$ satisfy
\begin{align*}
	\frac{2+p^{-1}}{pq-1}>\frac{n-1}{2},
\end{align*}
or it can be rewritten by
\begin{align*}
	\frac{q+2+p^{-1}}{pq-1}-\frac{n-1}{2}>\frac{q}{pq-1}=\alpha_{\mathrm{Shrift}}(p,q).
\end{align*}
The blow-up condition in the subcritical case for the weakly coupled system of wave equations with power nonlinearities is
\begin{align*}
	\frac{q+2+p^{-1}}{pq-1}-\frac{n-1}{2}>0.
\end{align*}
In other words, we may explain the blow-up condition for Nakao's type problem \eqref{Nakao's problem power} by a shrift curve with the magnitude $\alpha_{\mathrm{Shrift}}(p,q)$ of  the weakly coupled system of wave equations with power nonlinearities.

Indeed, this effect with the magnitude still holds for the Nakao's type problem with nonlinearities of derivative type. Let us turn to the blow-up result for the coupled system \eqref{Nakao's problem derivative}, i.e. $d_1=d_3=0$ and $d_2=d_4=1$ in the Cauchy problem \eqref{General Nakao's type problem}. Our main result in Theorem \ref{Thm 1} for $n\geqslant 2$ claims that the blow-up condition is
\begin{align*}
\frac{1}{pq-1}>\frac{n-1}{2},
\end{align*}
or it can be rewritten by
\begin{align*}
	\frac{q+1}{pq-1}-\frac{n-1}{2}>\frac{q}{pq-1}=\alpha_{\mathrm{Shrift}}(p,q).
\end{align*}
Therefore, we still may explain the blow-up condition for Nakao's type problem \eqref{Nakao's problem power} by a shrift magnitude $\alpha_{\mathrm{Shrift}}(p,q)$ of  the weakly coupled system of wave equations with nonlinearities of derivative type.

\section{Proof of Theorem \ref{Thm 1} via an iteration argument}\label{Proof of Thm 01}
\setcounter{equation}{0}

\subsection{Iteration frame}
In order to apply an iteration argument in the proof, we should derive integral inequalities for some suitable time-dependent functions. To begin with, let us introduce the eigenfunction $\Phi=\Phi(x)$ of the Laplace operator in $n$-dimensions Euclidean space such that
\begin{align*}
\Phi(x) & :=  \mathrm{e}^{x}+\mathrm{e}^{-x}  \qquad  \qquad \, \mbox{if} \ \ n=1, \\  \Phi(x) & := \int_{\mathbb{S}^{n-1}} \mathrm{e}^{x\cdot \omega}  \mathrm{d} \sigma_\omega \qquad \,\mbox{if} \ \  n\geqslant 2,
\end{align*}
where $\mathbb{S}^{n-1}$ is the $n-1$ dimensional sphere. This test function $\Phi$ has been introduced in the pioneering paper \cite{YordanovZhang2006}. It fulfills the property $\Delta\Phi=\Phi$ and the asymptotic behavior
\begin{align}\label{AsymptoticBehavior}
\Phi(x)\sim|x|^{-\frac{n-1}{2}}\mathrm{e}^{x}\ \ \mbox{as}\ \  |x|\to\infty.
\end{align}
Moreover, we define the test function $\Psi=\Psi(t,x)$ with separate variables such that $\Psi(t,x):= \mathrm{e}^{-t}\Phi(x)$. Clearly, the function $\Psi$ is a special solution of the homogeneous wave equation $\Psi_{tt}-\Delta\Psi=0$. By using asymptotic behavior \eqref{AsymptoticBehavior}, it immediately yields the estimate
\begin{align}\label{Est Psi}
\int_{B_{R+t}}\Psi(t,x)\mathrm{d}x\leqslant C_1(R+t)^{\frac{n-1}{2}}
\end{align}
for any $t\geqslant0$, where $C_1$ is a positive constant. The previous estimate \eqref{Est Psi} was shown in \cite{Lai-Takamura-2019}.

To construct the iteration frame, it is necessary for us to introduce some functionals with respect to $u_t$ and $v_t$ due to the derivative type nonlinearities of the hyperbolic coupled system \eqref{Nakao's problem derivative}. With the aid of the above test function $\Psi$, we may define new time-dependent functionals $F_1=F_1(t)$ and $F_2=F_2(t)$ such that
\begin{align*}
F_1(t):=\int_0^t\int_{\mb{R}^n}u_t(s,x)\Psi(s,x)\mathrm{d}x\mathrm{d}s\ \ \mbox{and} \ \ F_2(t):=\int_{\mb{R}^n}v_t(t,x)\Psi(t,x)\mathrm{d}x.
\end{align*}
Here, we should emphasize that $F'_1(t)$ has the similar form to $F_2(t)$, which is beneficial to process the iteration procedure later.

Due to the fact that $u,v$ are supported in a forward cone $\{(s,x)\in[0,t]\times\mb{R}^n:|x|\leqslant R+s\}$, we can apply the definition of energy solution $(u,v)$ with $\Psi$ to be the test function in \eqref{Defn.Energy.Solution 01} and \eqref{Defn.Energy.Solution 02}. For one thing, by using integration by parts in \eqref{Defn.Energy.Solution 01} with $\phi(t,x)=\Psi(t,x)$, we have
\begin{align*}
&\int_0^t\int_{\mb{R}^n}u_t(s,x)\Psi(s,x)\mathrm{d}x\mathrm{d}s+\int_{\mb{R}^n}(u_t(t,x)\Psi(t,x)+u(t,x)\Psi(t,x))\mathrm{d}x\\
&=\varepsilon\int_{\mb{R}^n}(u_0(x)+u_1(x))\Phi(x)\mathrm{d}x+\int_0^t\int_{\mb{R}^n}|v_t(s,x)|^p\Psi(s,x)\mathrm{d}x\mathrm{d}s,
\end{align*}
which can also be rewritten by
\begin{align}\label{Add 01}
&F_1'(t)+F_1(t)+\int_{\mb{R}^n}u(t,x)\Psi(t,x)\mathrm{d}x\notag\\
&=\varepsilon\int_{\mb{R}^n}(u_0(x)+u_1(x))\Phi(x)\mathrm{d}x+\int_0^t\int_{\mb{R}^n}|v_t(s,x)|^p\Psi(s,x)\mathrm{d}x\mathrm{d}s.
\end{align}
Taking time-derivative in the above equality and using $\Psi_t(t,x)=-\Psi(t,x)$ brings
\begin{align}\label{Add 02}
F_1''(t)+2F'_1(t)-\int_{\mb{R}^n}u(t,x)\Psi(t,x)\mathrm{d}x=\int_{\mb{R}^n}|v_t(t,x)|^p\Psi(t,x)\mathrm{d}x.
\end{align}
Adding up \eqref{Add 01} and \eqref{Add 02}, one may derive
\begin{align}\label{Eq F1 01}
F_1''(t)+3F_1'(t)+F_1(t)&=\varepsilon\int_{\mb{R}^n}(u_0(x)+u_1(x))\Phi(x)\mathrm{d}x+\int_0^t\int_{\mb{R}^n}|v_t(s,x)|^p\Psi(s,x)\mathrm{d}x\mathrm{d}s\notag\\
&\quad+\int_{\mb{R}^n}|v_t(t,x)|^p\Psi(t,x)\mathrm{d}x.
\end{align}
For another thing, we employ once integration by parts in \eqref{Defn.Energy.Solution 02} with $\psi(t,x)=\Psi(t,x)$ to get
%\begin{align*}
%\int_{\mb{R}^n}v_t(t,x)\Psi(t,x)\mathrm{d}x+\int_{\mb{R}^n}v(t,x)\Psi(t,x)\mathrm{d}x&=\varepsilon\int_{\mb{R}^n}(v_0(x)+v_1(x))\Phi(x)\mathrm{d}x\\
%&\quad+\int_{0}^t\int_{\mb{R}^n}|u_t(s,x)|^q\Psi(s,x)\mathrm{d}x\mathrm{d}s,
%\end{align*}
%which can be represented by
\begin{align*}
F_2(t)+\int_{\mb{R}^n}v(t,x)\Psi(t,x)\mathrm{d}x=\varepsilon\int_{\mb{R}^n}(v_0(x)+v_1(x))\Phi(x)\mathrm{d}x+\int_{0}^t\int_{\mb{R}^n}|u_t(s,x)|^q\Psi(s,x)\mathrm{d}x\mathrm{d}s.
\end{align*}
Similarly to the treatment of $F_1(t)$, we differentiate the last equality with respect to $t$, which implies
\begin{align*}
F'_2(t)+F_2(t)-\int_{\mb{R}^n}v(t,x)\Psi(t,x)\mathrm{d}x=\int_{\mb{R}^n}|u_t(t,x)|^q\Psi(t,x)\mathrm{d}x.
\end{align*}
Summarizing the derived equations, one has
\begin{align}\label{Eq F2 01}
F'_2(t)+2F_2(t)&=\varepsilon\int_{\mb{R}^n}(v_0(x)+v_1(x))\Phi(x)\mathrm{d}x+\int_{0}^t\int_{\mb{R}^n}|u_t(s,x)|^q\Psi(s,x)\mathrm{d}x\mathrm{d}s\notag\\
&\quad+\int_{\mb{R}^n}|u_t(t,x)|^q\Psi(t,x)\mathrm{d}x.
\end{align}

With the aim of constructing the iteration frame, we need to transfer \eqref{Eq F1 01} and \eqref{Eq F2 01} to suitable integral inequalities, respectively. Let us consider \eqref{Eq F1 01} initially. We now define a  time-dependent functional:
\begin{align*}
G_1(t):=F_1'(t)+\frac{3+\sqrt{5}}{2}F_1(t)-\varepsilon\int_{\mb{R}^n}u_1(x)\Phi(x)\mathrm{d}x-\int_0^t\int_{\mb{R}^n}|v_t(s,x)|^p\Psi(s,x)\mathrm{d}x\mathrm{d}s.
\end{align*}
Then, it is obvious from \eqref{Eq F1 01} that
\begin{align*}
G_1'(t)+\frac{3-\sqrt{5}}{2}G_1(t)&=F_1''(t)+3F_1'(t)+F_1(t)-\frac{(3-\sqrt{5})\varepsilon}{2}\int_{\mb{R}^n}u_1(x)\Phi(x)\mathrm{d}x\\
&\quad-\int_{\mb{R}^n}|v_t(t,x)|^p\Psi(t,x)\mathrm{d}x-\frac{3-\sqrt{5}}{2}\int_0^t\int_{\mb{R}^n}|v_t(s,x)|^p\Psi(s,x)\mathrm{d}x\mathrm{d}s\\
&=\varepsilon\int_{\mb{R}^n}u_0(x)\Phi(x)\mathrm{d}x+\frac{(\sqrt{5}-1)\varepsilon}{2}\int_{\mb{R}^n}u_1(x)\Phi(x)\mathrm{d}x\\
&\quad+\frac{\sqrt{5}-1}{2}\int_0^t\int_{\mb{R}^n}|v_t(s,x)|^p\Psi(s,x)\mathrm{d}x\mathrm{d}s.
\end{align*}
According to the nonnegative hypothesis on initial data $u_0$ and $u_1$, we are able to conclude 
\begin{align*}
\mathrm{e}^{-\frac{3-\sqrt{5}}{2}t}\left(\mathrm{e}^{\frac{3-\sqrt{5}}{2}t}G_1(t)\right)'=G_1'(t)+\frac{3-\sqrt{5}}{2}G_1(t)\geqslant0,
\end{align*}
which results
\begin{align*}
G_1(t)\geqslant\mathrm{e}^{-\frac{3-\sqrt{5}}{2}t}G_1(0)=\mathrm{e}^{-\frac{3-\sqrt{5}}{2}t}\left(F_1'(0)+\frac{3+\sqrt{5}}{2}F_1(0)-\varepsilon\int_{\mb{R}^n}u_1(x)\Phi(x)\mathrm{d}x\right)=0.
\end{align*}
For this reason, we obtain
\begin{align}\label{Sup 01}
\mathrm{e}^{-\frac{3+\sqrt{5}}{2}t}\left(\mathrm{e}^{\frac{3+\sqrt{5}}{2}t}F_1(t)\right)'\geqslant\varepsilon\int_{\mb{R}^n}u_1(x)\Phi(x)\mathrm{d}x+\int_0^t\int_{\mb{R}^n}|v_t(s,x)|^p\Psi(s,x)\mathrm{d}x\mathrm{d}s.
\end{align}
By ignoring the nonnegative nonlinear integral term on the right-hand side, multiplying the previous equality by $\mathrm{e}^{\frac{3+\sqrt{5}}{2}t}$ and integrating the resultant over $[0,t]$, we arrive at
\begin{align}\label{Pre First Lower F1 01}
F_1(t)\geqslant\mathrm{e}^{-\frac{3+\sqrt{5}}{2}t}F_1(0)+\frac{2\varepsilon}{3+\sqrt{5}}\left(1-\mathrm{e}^{-\frac{3+\sqrt{5}}{2}t}\right)\int_{\mb{R}^n}u_1(x)\Phi(x)\mathrm{d}x\geqslant C_2\varepsilon>0
\end{align}
for any $t\geqslant1$, where $C_2$ is a suitably positive constant depending on $u_1$. Here, we used nontrivial assumption on $u_1$ and $F_1(0)=0$. What's more, by omitting the term containing initial data $u_1$ in \eqref{Sup 01} we find that
\begin{align}\label{Iteration Frame 01}
F_1(t)&\geqslant\int_0^t\mathrm{e}^{\frac{3+\sqrt{5}}{2}(\tau-t)}\int_0^{\tau}\int_{\mb{R}^n}|v_t(s,x)|^p\Psi(s,x)\mathrm{d}x\mathrm{d}s\mathrm{d}\tau\notag\\
&\geqslant\int_0^t\mathrm{e}^{\frac{3+\sqrt{5}}{2}(\tau-t)}\int_0^{\tau}|F_2(s)|^p\left(\int_{|x|\leqslant R+s}\Psi(s,x)\mathrm{d}x\right)^{-(p-1)}\mathrm{d}s\mathrm{d}\tau\notag\\
&\geqslant C_1^{1-p}\int_0^t\mathrm{e}^{\frac{3+\sqrt{5}}{2}(\tau-t)}\int_0^{\tau}(R+s)^{-\frac{(n-1)(p-1)}{2}}|F_2(s)|^p\mathrm{d}s\mathrm{d}\tau,
\end{align}
where we used the support condition for the wave model and the estimate \eqref{Est Psi}.\\
Next, we treat \eqref{Eq F2 01} by constructing another time-dependent functional such that
\begin{align}\label{Defn G2}
G_2(t):=F_2(t)-\frac{\varepsilon}{2}\int_{\mb{R}^n}v_1(x)\Phi(x)\mathrm{d}x-\frac{1}{2}\int_0^t\int_{\mb{R}^n}|u_t(s,x)|^q\Psi(s,x)\mathrm{d}x\mathrm{d}s.
\end{align}
In other words, in the light of \eqref{Eq F2 01} we find
\begin{align*}
\mathrm{e}^{-2t}\left(\mathrm{e}^{2t}G_2(t)\right)'&=F_2'(t)+2F_2(t)-\varepsilon\int_{\mb{R}^n}v_1(x)\Phi(x)\mathrm{d}x\\
&\quad-\frac{1}{2}\int_{\mb{R}^n}|u_t(t,x)|^q\Psi(t,x)\mathrm{d}x-\int_0^t\int_{\mb{R}^n}|v_t(s,x)|^q\Psi(s,x)\mathrm{d}x\mathrm{d}s\\
&=\varepsilon\int_{\mb{R}^n}v_0(x)\Phi(x)\mathrm{d}x+\frac{1}{2}\int_{\mb{R}^n}|u_t(t,x)|^q\Psi(t,x)\mathrm{d}x\geqslant0,
\end{align*}
where the nonnegativity of $v_0$ was applied. It immediately conduces to
\begin{align*}
G_2(t)\geqslant \mathrm{e}^{-2t}G_2(0)=\frac{\mathrm{e}^{-2t}\varepsilon}{2}\int_{\mb{R}^n}v_1(x)\Phi(x)\mathrm{d}x\geqslant0
\end{align*}
from the nonnegativity of $v_1$.  Consequently, the nontrivial assumption on $v_1$ associated with the relation \eqref{Defn G2} shows
 \begin{align}\label{Pre First Lower F1 02}
 F_2(t)\geqslant\frac{\varepsilon}{2}\int_{\mb{R}^n}v_1(x)\Phi(x)\mathrm{d}x=C_3\varepsilon>0,
 \end{align}
 with a positive constant $C_3$ depending on $v_1$, and
\begin{align}\label{Iteration Frame 02}
F_2(t)&\geqslant\frac{1}{2}\int_0^t\int_{\mb{R}^n}|u_t(s,x)|^q\Psi(s,x)\mathrm{d}x\mathrm{d}s\notag\\
&\geqslant\frac{C_1^{1-q}}{2}\int_0^t(R+s)^{-\frac{(n-1)(q-1)}{2}}|F_1'(s)|^q\mathrm{d}s\notag\\
&\geqslant\frac{C_1^{1-q}}{2}\left(\int_0^t(R+s)^{\frac{n-1}{2}}\mathrm{d}s\right)^{-(q-1)}\left(\int_0^t|F_1'(s)|\mathrm{d}s\right)^q\notag\\
&\geqslant C_4(R+t)^{-\frac{(n+1)(q-1)}{2}}|F_1(t)|^q,
\end{align}
with a positive constant $C_4$, where we utilized H\"older's inequality 
\begin{align*}
\int_0^t|F_1'(s)|\mathrm{d}s\leqslant\left(\int_0^t|F_1'(s)|^q(R+s)^{-\frac{(n-1)(q-1)}{2}}\mathrm{d}s\right)^{\frac{1}{q}}\left(\int_0^t(R+s)^{\frac{(n-1)(q-1)q'}{2q}}\mathrm{d}s\right)^{\frac{1}{q'}},
\end{align*}
associated with \eqref{Est Psi} again and
\begin{align*}
|F_1(t)|=|F_1(t)-F_1(0)|=\left|\int_0^tF_1'(s)\mathrm{d}s\right|\leqslant\int_0^t|F'_1(s)|\mathrm{d}s.
\end{align*}

The further step is to investigate first lower bound estimates for the functionals $F_1(t)$ as well as $F_2(t)$, individually. On one hand, we combine \eqref{Pre First Lower F1 02} with \eqref{Iteration Frame 01} to deduce
\begin{align*}
F_1(t)&\geqslant C_1^{1-p}C_3^p\varepsilon^p\int_0^t\mathrm{e}^{\frac{3+\sqrt{5}}{2}(\tau-t)}\int_0^{\tau}(R+s)^{-\frac{(n-1)(p-1)}{2}}\mathrm{d}s\mathrm{d}\tau\\
&\geqslant C_1^{1-p}C_3^p\varepsilon^p(R+t)^{-\frac{(n-1)(p-1)}{2}}\int_{t/2}^t\mathrm{e}^{\frac{3+\sqrt{5}}{2}(\tau-t)}\tau\mathrm{d}\tau\\
&\geqslant\frac{C_1^{1-p}C_3^p}{3+\sqrt{5}}\varepsilon^p(R+t)^{-\frac{(n-1)(p-1)}{2}}t\left(1-\mathrm{e}^{-\frac{3+\sqrt{5}}{4}t}\right)\\
&\geqslant\frac{C_1^{1-p}C_3^p}{3+\sqrt{5}}\left(1-\mathrm{e}^{-1-\frac{3+\sqrt{5}}{4}}\right)\varepsilon^p(R+t)^{-\frac{(n-1)(p-1)}{2}}(t-L_1)
\end{align*}
for any $t\geqslant L_1:=1+4/(3+\sqrt{5})$, which provides first lower bound estimates for $F_1(t)$. On the other hand, we summarize \eqref{Pre First Lower F1 01} and \eqref{Iteration Frame 02}. It results
\begin{align*}
F_2(t)\geqslant C_2^qC_4\varepsilon^q(R+t)^{-\frac{(n+1)(q-1)}{2}}
\end{align*}
for any $t\geqslant L_1>1$. This choice of $L_1$ is concerned about the slicing procedure dealing with the unbounded multiplier in the next subsection.

All in all, we derived first lower bound estimates as follows:
\begin{align}
F_1(t)&\geqslant D_1(R+t)^{-\alpha_1}(t-L_1)^{\beta_1},\label{First Lower B F1}\\
F_2(t)&\geqslant Q_1(R+t)^{-a_1}(t-L_1)^{b_1},\label{First Lower B F2}
\end{align}
for any $t\geqslant L_1$, where the multiplicative constants are given by
\begin{align*}
D_1:=\frac{C_1^{1-p}C_3^p}{3+\sqrt{5}}\left(1-\mathrm{e}^{-1-\frac{3+\sqrt{5}}{4}}\right)\varepsilon^p, \ \ Q_1:=C_2^qC_4\varepsilon^q,
\end{align*}
and the exponents are represented by
\begin{align*}
\alpha_1:=\frac{(n-1)(p-1)}{2}, \ \ a_1:=\frac{(n+1)(q-1)}{2},\ \ \beta_1:=1,\ \ b_1:=0.
\end{align*}
We remark that the above constants are nonnegative.

\subsection{Iteration argument}
In this part, we will derive sequences of lower bound estimates for the functionals $F_1(t)$ and $F_2(t)$ by using some derived inequalities in the last subsection. To be specific, the following lower bounds will be proved:
\begin{align}
F_1(t)&\geqslant D_j(R+t)^{-\alpha_j}(t-L_j)^{\beta_j},\label{Seq F1}\\
F_2(t)&\geqslant Q_j(R+t)^{-a_j}(t-L_j)^{b_j},\label{Seq F2}
\end{align}
for any $t\geqslant L_j$, where $\{D_j\}_{j\geqslant 1}$, $\{Q_j\}_{j\geqslant 1}$, $\{\alpha_j\}_{j\geqslant 1}$, $\{a_j\}_{j\geqslant 1}$, $\{\beta_j\}_{j\geqslant 1}$ and $\{b_j\}_{j\geqslant 1}$ are sequences of nonnegative real numbers that will be determined later in the iteration procedure. Motivated by the recent papers \cite{ChenPalmieri201901,ChenPalmieri201902}, we may define a crucial sequence $\{L_j\}_{j\geqslant 1}$ of the partial products of the convergent infinite product
\begin{align}\label{Seq_Nakao}
\prod\limits_{k=1}^{\infty}\ell_{k}\ \ \mbox{with}\ \ \ell_k:=1+\frac{4}{3+\sqrt{5}}(pq)^{-\frac{k-1}{2}}\ \ \mbox{for any}\ \  k\geqslant1,
\end{align}
that is,
\begin{align}\label{Seq_Nakao_2}
L_j:=\prod\limits_{k=1}^{j}\ell_{k}\ \ \mbox{for any}\ \  j\geqslant1.
\end{align}
Here, we recall that $L_1=\ell_1=1+4/(3+\sqrt{5})$. Essentially, thanks to the ratio test and
\begin{align*}
\lim\limits_{k\to\infty}\frac{\ln\ell_{k+1}}{\ln \ell_k}=\lim\limits_{k\to\infty}\frac{\ell_k}{\ell_{k+1}(pq)^{1/2}}=(pq)^{-1/2}<1,
\end{align*}
we claim that the infinite product
\begin{align*}
\prod\limits_{k=1}^{\infty}\ell_{k}=\exp\left(\sum\limits_{k=1}^{\infty}\ln\ell_k\right)
\end{align*}
 is convergent. Furthermore, the desired estimates \eqref{Seq F1} and \eqref{Seq F2} for $j=1$ are given in \eqref{First Lower B F1} and \eqref{First Lower B F2}, respectively.
 
 As a consequence, with the aim of demonstrating \eqref{Seq F1} and \eqref{Seq F2}, we just need to procure the induction step with the aim of proving \eqref{Seq F1} and \eqref{Seq F2}. In other words, by assuming that \eqref{Seq F1} and \eqref{Seq F2} hold for $j$, one oughts to prove them being valid for $j+1$. Let us first substitute \eqref{Seq F2} into \eqref{Iteration Frame 01}, which leads to
 \begin{align*}
 F_1(t)&\geqslant C_1^{1-p}Q_j^p\int_0^t\mathrm{e}^{\frac{3+\sqrt{5}}{2}(\tau-t)}\int_0^{\tau}(R+s)^{-\frac{(n-1)(p-1)}{2}-a_jp}(s-L_j)^{b_jp}\mathrm{d}s\mathrm{d}\tau\\
 &\geqslant C_1^{1-p}Q_j^p(R+t)^{-\frac{(n-1)(p-1)}{2}-a_jp}\int_{L_j}^t\mathrm{e}^{\frac{3+\sqrt{5}}{2}(\tau-t)}\int_{L_j}^{\tau}(s-L_j)^{b_jp}\mathrm{d}s\mathrm{d}\tau\\
 &\geqslant \frac{C_1^{1-p}Q_j^p}{b_jp+1}(R+t)^{-\frac{(n-1)(p-1)}{2}-a_jp}\int_{L_j}^t\mathrm{e}^{\frac{3+\sqrt{5}}{2}(\tau-t)}(\tau-L_j)^{b_jp+1}\mathrm{d}\tau.
 \end{align*}
 In view of $t\geqslant L_{j+1}=L_j\ell_{j+1}$, i.e. $L_j\leqslant t/\ell_{j+1}$, we may instantly shrink the interval $[L_j,t]$ into $[t/\ell_{j+1},t]$ so that
 \begin{align*}
 F_1(t)&\geqslant \frac{C_1^{1-p}Q_j^p}{b_jp+1}(R+t)^{-\frac{(n-1)(p-1)}{2}-a_jp}\int_{t/\ell_{j+1}}^t\mathrm{e}^{\frac{3+\sqrt{5}}{2}(\tau-t)}(\tau-L_j)^{b_jp+1}\mathrm{d}\tau\\
 &\geqslant\frac{2C_1^{1-p}Q_j^p}{(3+\sqrt{5})(b_jp+1)\ell_{j+1}^{b_jp+1}}\left(1-\mathrm{e}^{\frac{3+\sqrt{5}}{2}(1/\ell_{j+1}-1)t}\right)(R+t)^{-\frac{(n-1)(p-1)}{2}-a_jp}(t-L_{j+1})^{b_jp+1}.
 \end{align*}
 By considering $t\geqslant L_{j+1}\geqslant\ell_{j+1}$ with the formula of $\ell_{j+1}$, one observes
 \begin{align*}
 1-\mathrm{e}^{\frac{3+\sqrt{5}}{2}(1/\ell_{j+1}-1)t}&\geqslant 1-\mathrm{e}^{-\frac{3+\sqrt{5}}{2}(\ell_{j+1}-1)}\geqslant\frac{3+\sqrt{5}}{2}(\ell_{j+1}-1)\left(1-\frac{3+\sqrt{5}}{4}(\ell_{j+1}-1)\right)\\
 &\geqslant 2(pq)^{-\frac{j}{2}}\left(1-(pq)^{-\frac{j}{2}}\right)\geqslant 2\left((pq)^{\frac{1}{2}}-1\right)(pq)^{-j}>0
 \end{align*}
for any $j\geqslant 1$. In conclusion, it yields
 \begin{align*}
 F_1(t)\geqslant\frac{4C_1^{1-p}\left((pq)^{\frac{1}{2}}-1\right)(pq)^{-j}Q_j^p}{(3+\sqrt{5})(b_jp+1)\ell_{j+1}^{b_jp+1}}(R+t)^{-\frac{(n-1)(p-1)}{2}-a_jp}(t-L_{j+1})^{b_jp+1}
 \end{align*}
 for any $t\geqslant L_{j+1}$. Then, the combination of \eqref{Iteration Frame 02} as well as \eqref{Seq F1} shows
 \begin{align*}
 F_2(t)\geqslant C_4D_j^q(R+t)^{-\frac{(n+1)(q-1)}{2}-\alpha_jq}(t-L_{j+1})^{\beta_jq}
 \end{align*}
 for any $t\geqslant L_{j+1}$, where we used the fact that $L_j\leqslant L_j\ell_{j+1}=L_{j+1}$ with $\ell_{j+1}>1$.
 
 In other words, \eqref{Seq F1} and \eqref{Seq F2} are valid supposing that
 \begin{align*}
 D_{j+1}&:=\frac{4C_1^{1-p}\left((pq)^{\frac{1}{2}}-1\right)(pq)^{-j}}{(3+\sqrt{5})(b_jp+1)\ell_{j+1}^{b_jp+1}}Q_j^p, \ \ Q_{j+1}:=C_4D_j^q,\\
 \alpha_{j+1}:=\frac{(n-1)(p-1)}{2}+&a_jp,\ \ a_{j+1}:=\frac{(n+1)(q-1)}{2}+\alpha_jq,\ \ \beta_{j+1}=1+b_jp, \ \ b_{j+1}:=\beta_jq.
 \end{align*}
\subsection{Upper bound estimates for the lifespan}
In the last subsection, we derive a sequence of lower bound estimates for $F_1(t)$ and $F_2(t)$, respectively. In the forthcoming part, we will demonstrate that the $j$-dependent lower bounds for the functionals $F_1(t)$ and $F_2(t)$ blow up as $j\to\infty$. At the same time, the blow-up result and upper bound estimates for the lifespan stated in Theorem \ref{Thm 1} will be concluded.

We will begin with the explicit formulas for the sequences $\alpha_j,\beta_j,a_j,b_j$, which devote to estimates for the multiplicative constants $D_j$ and $Q_j$. 

Particularly, concerning the formulas of $\alpha_j$ and $a_j$, we need to discuss the case when $j$ is an odd integer only, which is sufficient for our proof. Taking account of the relation between $\alpha_j$ and $a_j$, we may get for odd number $j$ that
\begin{align*}
\alpha_j&=\frac{(n-1)(p-1)}{2}+a_{j-1}p=\frac{(n+1)pq-2p-(n-1)}{2}+\alpha_{j-2}pq\\
&=\frac{(n+1)pq-2p-(n-1)}{2}\sum\limits_{k=0}^{(j-3)/2}(pq)^k+\alpha_1(pq)^{\frac{j-1}{2}}\\
&=\left(\alpha_1+\frac{(n+1)pq-2p-(n-1)}{2(pq-1)}\right)(pq)^{\frac{j-1}{2}}-\frac{(n+1)pq-2p-(n-1)}{2(pq-1)},
\end{align*}
and similarly,
\begin{align*}
a_j&=\frac{(n+1)(q-1)}{2}+\alpha_{j-1}q=\frac{(n-1)pq+2q-(n+1)}{2}+a_{j-2}pq\\
&=\frac{(n-1)pq+2q-(n+1)}{2}\sum\limits_{k=0}^{(j-3)/2}(pq)^k+a_1(pq)^{\frac{j-1}{2}}\\
&=\left(a_1+\frac{(n-1)pq+2q-(n+1)}{2(pq-1)}\right)(pq)^{\frac{j-1}{2}}-\frac{(n-1)pq+2q-(n+1)}{2(pq-1)}.
\end{align*}
Furthermore, by the definition of $\beta_j$ and $b_j$, one derives for odd number $j$ that
\begin{align*}
\beta_j&=1+b_{j-1}p=1+\beta_{j-2}pq=\sum\limits_{k=0}^{(j-3)/2}(pq)^k+\beta_1(pq)^{\frac{j-1}{2}}=\left(\beta_1+\frac{1}{pq-1}\right)(pq)^{\frac{j-1}{2}}-\frac{1}{pq-1},\\
b_j&=\beta_{j-1}q=q+b_{j-2}pq=q\sum\limits_{k=0}^{(j-3)/2}(pq)^k+b_1(pq)^{\frac{j-1}{2}}=\left(b_1+\frac{q}{pq-1}\right)(pq)^{\frac{j-1}{2}}-\frac{q}{pq-1}.
\end{align*}
For an even number $j$, which means that $j-1$ is an odd number, we make use of the previous two equalities to arrive at
\begin{align*}
\beta_j&=1+b_{j-1}p=q^{-1}\left(b_1+\frac{q}{pq-1}\right)(pq)^{\frac{j}{2}}-\frac{1}{pq-1},\\
b_j&=\beta_{j-1}q=p^{-1}\left(\beta_1+\frac{1}{pq-1}\right)(pq)^{\frac{j}{2}}-\frac{q}{pq-1}.
\end{align*}
For this reason, it holds
\begin{align*}
\beta_j\leqslant B_0(pq)^{\frac{j}{2}}\ \ \mbox{and}\ \ b_j\leqslant B_1(pq)^{\frac{j}{2}}
\end{align*}
for any $j\geqslant1$, where $B_0=B_0(p,q,n)$ and $B_1=B_1(p,q,n)$ are positive constants independent of $j$.

Before estimating the constants $D_j$ and $Q_j$ from the below, we apply  L'H\^opital's rule to show
\begin{align*}
\lim\limits_{j\to\infty}\ell_j^{b_{j-1}p+1}&=\lim\limits_{j\to\infty}\ell_j^{\beta_j}\leqslant\lim\limits_{j\to\infty}\exp\left(B_0(pq)^{\frac{j}{2}}\ln\left(1+\frac{4}{3+\sqrt{5}}(pq)^{-\frac{j-1}{2}}\right)\right)\\
&=\exp\left(\frac{4B_0}{3+\sqrt{5}}(pq)^{\frac{1}{2}}\right)>0
\end{align*}
so that there exists a suitable constant satisfying $1/\ell_j^{b_{j-1}p+1}\geqslant M>0$ for any $j\geqslant 1$. As a result, the next iterated relations for the lower bounds come:
\begin{align}\label{Iter Dj}
D_j&=\frac{4C_1^{1-p}\left((pq)^{1/2}-1\right)(pq)^{-j+1}}{(3+\sqrt{5})(b_{j-1}p+1)\ell_j^{b_{j-1}p+1}}Q_{j-1}^p\geqslant\frac{4C_1^{1-p}\left((pq)^{1/2}-1\right)M}{(3+\sqrt{5})B_0}(pq)^{-\frac{3}{2}j+1}Q_{j-1}^p\notag\\
&\geqslant\frac{4C_1^{1-p}C_4^p\left((pq)^{1/2}-1\right)M}{(3+\sqrt{5})B_0}(pq)^{-\frac{3}{2}j+1}D_{j-2}^{pq}=:E_0(pq)^{-\frac{3}{2}j+1}D_{j-2}^{pq},
\end{align}
and simultaneously,
\begin{align}\label{Iter Qj}
Q_j&=C_4D_{j-1}^q\geqslant\frac{4^qC_1^{(1-p)q}C_4\left((pq)^{1/2}-1\right)^qM^q}{(3+\sqrt{5})^qB_0^q}(pq)^{-\frac{3}{2}jq+q}Q_{j-2}^{pq}=: E_1(pq)^{-\frac{3}{2}jq+q}Q_{j-2}^{pq},
\end{align}
with suitable constants $E_0=E_0(p,q,n)>0$ and $E_1=E_1(p,q,n)>0$ independent of $j$.\\
 Considering \eqref{Iter Dj} with an odd number $j$, we take the logarithmic on the both sides to deduce
\begin{align*}
\log D_j&\geqslant pq\log D_{j-2}-\left(\frac{3}{2}j-1\right)\log (pq)+\log E_0\\
&\geqslant (pq)^{\frac{j-1}{2}}\log D_1-\frac{3}{2}\log(pq)\sum\limits_{k=1}^{(j-1)/2}\left((j+2-2k)(pq)^{k-1}\right)\\
&\quad+\left(\log(pq)+\log E_0\right)\sum\limits_{k=1}^{(j-1)/2}(pq)^{k-1}\\
&=(pq)^{\frac{j-1}{2}}\left(\log D_1+\frac{\log (pq)}{2(pq-1)^2}(1-7pq)+\frac{\log E_0}{pq-1}\right)\\
&\quad+\frac{\log(pq)}{pq-1}\left(\frac{3}{2}\left(\frac{2pq}{pq-1}+j\right)-1\right)-\frac{\log E_0}{pq-1},
\end{align*}
where we used the next formula in the last line of the chain inequalities:
\begin{align*}
\sum\limits_{k=1}^{(j-1)/2}\left((j+2-2k)(pq)^{k-1}\right)=\frac{1}{pq-1}\left(\frac{2pq}{pq-1}\left(\frac{3}{2}(pq)^{\frac{j-1}{2}}-\frac{1}{2}(pq)^{\frac{j-3}{2}}-1\right)-j\right).
\end{align*}
Thus, for all nonnegative odd numbers $j$ satisfying
\begin{align*}
j\geqslant j_0:=\left\lceil\frac{2}{3}+\frac{2\log E_0}{3\log(pq)}-\frac{2pq}{pq-1}\right\rceil,
\end{align*}
we conclude
\begin{align*}
\log D_j&\geqslant (pq)^{\frac{j-1}{2}}\left(\log D_1+\frac{\log (pq)}{2(pq-1)^2}(1-7pq)+\frac{\log E_0}{pq-1}\right)\\
&= (pq)^{\frac{j-1}{2}}\log\left(D_1(pq)^{(1-7pq)/(2(pq-1)^2)}E_0^{1/(pq-1)}\right)=(pq)^{\frac{j-1}{2}}\log(E_2\varepsilon^p)
\end{align*}
for a suitable positive constant $E_2=E_2(p,q,n)$. By the same way of calculation, we may illustrate
\begin{align*}
\log Q_j&\geqslant pq\log Q_{j-2}-\left(\frac{3}{2}jq-q\right)\log(pq)+\log E_1\\
&\geqslant(pq)^{\frac{j-1}{2}}\log Q_1-\frac{3}{2}q\log(pq)\sum\limits_{k=1}^{(j-1)/2}\left((j+2-2k)(pq)^{k-1}\right)\\
&\quad+\left(q\log(pq)+\log E_1\right)\sum\limits_{k=1}^{(j-1)/2}(pq)^{k-1}\\
&=(pq)^{\frac{j-1}{2}}\left(\log Q_1+\frac{q\log(pq)}{2(pq-1)^2}(1-7pq)+\frac{\log E_1}{pq-1}\right)\\
&\quad+\frac{\log(pq)}{pq-1}\left(\frac{3q}{2}\left(\frac{2pq}{pq-1}+j\right)-q\right)-\frac{\log E_1}{pq-1}.
\end{align*}
Consequently, for all nonnegative odd numbers $j$ fulfilling
\begin{align*}
j\geqslant j_1:=\left\lceil\frac{2}{3}+\frac{2\log E_1}{3q\log(pq)}-\frac{2pq}{pq-1}\right\rceil,
\end{align*}
we conclude
\begin{align*}
\log Q_j&\geqslant (pq)^{\frac{j-1}{2}}\left(\log Q_1+\frac{q\log(pq)}{2(pq-1)^2}(1-7pq)+\frac{\log E_1}{pq-1}\right)\\
&=(pq)^{\frac{j-1}{2}}\log\left(Q_1(pq)^{q(1-7pq)/(2(pq-1)^2)}E_1^{1/(pq-1)}\right)=(pq)^{\frac{j-1}{2}}\log(E_3\varepsilon^q)
\end{align*}
for a suitable positive constant $E_3=E_3(p,q,n)$.

Let us now denote
\begin{align*}
L:= \lim\limits_{j\to\infty}L_j=\prod\limits_{j=1}^{\infty}\ell_j>1.
\end{align*}
Due to $\ell_j>1$ the sequence $\{L_j\}_{j\geqslant1}$ is converging to $L$ as $j\to\infty$. Namely, the relations \eqref{Seq F1} and \eqref{Seq F2} hold for any odd number $j\geqslant1$ and any $t\geqslant L$.

Let us now consider an odd number $j$ such that $j\geqslant \max\{j_0,j_1\}$. The estimate \eqref{Seq F1} can be shown by
\begin{align*}
F_1(t)&\geqslant\exp\left((pq)^{\frac{j-1}{2}}\log(E_2\varepsilon^p)\right)(R+t)^{-\alpha_j}(t-L)^{\beta_j}\\
&\geqslant\exp\left((pq)^{\frac{j-1}{2}}\left(\log(E_2\varepsilon^p)-\left(\alpha_1+\tfrac{(n+1)pq-2p-(n-1)}{2(pq-1)}\right)\log(R+t)+\left(\beta_1+\tfrac{1}{pq-1}\right)\log(t-L)\right)\right)\\
&\quad\times(R+t)^{\frac{(n+1)pq-2p-(n-1)}{2(pq-1)}}(t-L)^{-\frac{1}{pq-1}}
\end{align*}
for any odd number $j\geqslant\max\{j_0,j_1\}$ and any $t\geqslant L$.  Choosing $t\geqslant\max \{R,2L\}$, since $R+t\leqslant 2t$ and $t-L\geqslant t/2$, the functional $F_1(t)$ can be estimated by the following way:
\begin{align}
F_1(t)&\geqslant\exp\left((pq)^{\frac{j-1}{2}}\log\left(E_2\varepsilon^p2^{-\frac{(n-1)(p-1)}{2}-\frac{(n+1)pq-2p-(n-1)}{2(pq-1)}-\frac{pq}{pq-1}}t^{-\frac{(n-1)(p-1)}{2}-\frac{(n+1)pq-2p-(n-1)}{2(pq-1)}+\frac{pq}{pq-1}}\right)\right)\notag\\
&\quad\times(R+t)^{\frac{(n+1)pq-2p-(n-1)}{2(pq-1)}}(t-L)^{-\frac{1}{pq-1}}\label{L.B.U(t)}
\end{align}
for any odd number $j\geqslant\max\{j_0,j_1\}$. The exponent of $t$, in the previous one, can be represented as follows:
\begin{align*}
&-\frac{(n-1)(p-1)}{2}-\frac{(n+1)pq-2p-(n-1)}{2(pq-1)}+\frac{pq}{pq-1}\\
&=\frac{p(-(n-1)pq+n+1)}{2(pq-1)}=:pT_1(p,q,n).
\end{align*}
By our assumption that $pq<(n+1)/(n-1)$ for any $n\geqslant 2$ and $p,q>1$ for any $n\geqslant1$, the power of $t$ in the exponential term of \eqref{L.B.U(t)} is positive.

In a similar way to the above, we may deduce the lower bound estimate for an odd number $j$ fulfilling $j\geqslant \max\{j_0,j_1\}$
\begin{align}
F_2(t)&\geqslant\exp\left((pq)^{\frac{j-1}{2}}\log\left(E_3\varepsilon^q2^{-\frac{(n+1)(q-1)}{2}-\frac{(n-1)pq+2q-(n+1)}{2(pq-1)}-\frac{q}{pq-1}}t^{-\frac{(n+1)(q-1)}{2}-\frac{(n-1)pq+2q-(n+1)}{2(pq-1)}+\frac{q}{pq-1}}\right)\right)\notag\\
&\quad\times(R+t)^{\frac{(n-1)pq+2q-(n+1)}{2(pq-1)}}(t-L)^{-\frac{q}{pq-1}}.\label{L.B.V(t)}
\end{align}
Thus, the power of $t$ in the exponential term can be represented by
\begin{align*}
&-\frac{(n+1)(q-1)}{2}-\frac{(n-1)pq+2q-(n+1)}{2(pq-1)}+\frac{q}{pq-1}\\
&=\frac{q(-(n+1)pq+2p+n+1)}{2(pq-1)}=:qT_2(p,q,n).
\end{align*}
 By assuming $(n+1)pq-2p-(n+1)<0$, the power for $t$ in the exponential term of \eqref{L.B.V(t)} is positive. We should emphasize that
 \begin{align*}
 \left\{(p,q):(n+1)pq-2p-(n+1)<0\right\}\subseteq\left\{(p,q):(n-1)pq<(n+1)\right\}
 \end{align*}
for all $p,q>1$ and $n\geqslant 1$. To put it differently, the condition $(n-1)pq<n+1$ is sufficient to guarantee the positivity of the power for $t$ in the exponential term of \eqref{L.B.V(t)}.

Eventually, for studying upper bound estimates for the lifespan, we now should introduce $\varepsilon_0=\varepsilon_0(u_0,u_1,v_0,v_1,p,q,n,R)>0$ such that
\begin{align*}
\left(E_2^{-1}2^{\frac{(n-1)(p-1)}{2}+\frac{(n+1)pq-2p-(n-1)}{2(pq-1)}+\frac{pq}{pq-1}}\right)^{1/(pT_1(p,q,n))}:= E_4\geqslant \varepsilon_0^{1/T_1(p,q,n)}.
\end{align*}
Hence, for $\varepsilon\in(0,\varepsilon_0]$ and for  $t>E_3\varepsilon^{-1/T_1(p,q,n)}$ as well as $t\geqslant \max\{R,2L\}$, letting $j\to\infty$ in \eqref{L.B.U(t)}, we claim that the lower bound for the functional $F_1(t)$ blows up. By the same way, in the case when $T_2(p,q,n)>0$, then we also can find a positive constant $\varepsilon_0=\varepsilon_0(u_0,u_1,v_0,v_1,p,q,n,R)>0$ such that
\begin{align*}
\left(E_3^{-1}2^{\frac{(n+1)(q-1)}{2}+\frac{(n-1)pq+2q-(n+1)}{2(pq-1)}+\frac{q}{pq-1}}\right)^{1/(qT_2(p,q,n))}:= E_5\geqslant \varepsilon_0^{1/T_2(p,q,n)}.
\end{align*}
For $\varepsilon\in(0,\varepsilon_0]$ and $t>E_5\varepsilon^{-1/T_2(p,q,n)}$ carrying $t\geqslant \max\{R,2L\}$, letting $j\to\infty$ in \eqref{L.B.V(t)}, we may immediately show that the lower bound for the functional $F_2(t)$ blows up. In conclusion, these statements proved that the energy solution $(u,v)$ is not defined globally in time and, simultaneously, the lifespan of this local (in time) solution $(u,v)$ can be estimated by
\begin{align*}
T(\varepsilon)\leqslant C\varepsilon^{-1/\max\{T_1(p,q,n),T_2(p,q,n)\}}=C\varepsilon^{-1/T_1(p,q,n)},
\end{align*}
where we used $T_1(p,q,n)> T_2(p,q,n)$ for any $n\geqslant1$. The proof of the theorem is complete.

\section*{Acknowledgments}
The author thanks  Michael Reissig (TU Bergakademie Freiberg) for the suggestions in the preparation of the paper.
%\bibliography{References}

% ------------------------------------------------------------------------
\end{document}